\pdfoutput=1
\documentclass[11pt]{amsart}
\usepackage{amsmath,amssymb,amsthm}
\usepackage[margin=1in]{geometry}
\usepackage{hyperref}
\usepackage{booktabs}

\newtheorem{theorem}{Theorem}
\newtheorem{lemma}[theorem]{Lemma}

\newtheorem{corollary}[theorem]{Corollary}
\theoremstyle{definition}

\theoremstyle{remark}
\newtheorem{remark}[theorem]{Remark}

\newcommand{\fl}[1]{\left\lfloor #1\right\rfloor}
\newcommand{\rep}{\operatorname{rep}}
\newcommand{\sh}{\operatorname{sh}}
\hypersetup{pdftitle={Kimberling's linear complementary equations: automatic structure and exact bounds}, pdfauthor={Alex Ashburn}}

\title[Kimberling's linear complementary equations]{Kimberling's linear complementary equations:\\ automatic structure and exact bounds}

\author{Alex Ashburn}
\address{Independent Researcher}
\email{aashburn1998@gmail.com}


\begin{document}

\begin{abstract}
In 2018 Kimberling added to the On-Line Encyclopedia of Integer Sequences the solutions of complementary equations such as $a(n)=2b(n-1)-b(n-2)+2n$, where $a(0)=1$, $a(1)=2$, $b(0)=3$, $b(1)=4$, and $b$ is the increasing sequence of positive integers not in $a$.
For $21$ of these sequences he conjectured bounds of the form $L<x(n)-\rho n<U$ for $n\ge1$, where $\rho$ is the growth constant of $x$.
The $11$ equations behind these conjectures have the form $a(n)=2b(n-s)-b(n-s-1)+kn+c$ with $s\in\{0,1\}$ and $k\in\{2,3,4\}$.
For each of them we show that $a(n)-(k+1)n$ and $b(n)-n$ are sums of an explicit shift function and an automatic sequence in the Ostrowski numeration system of $[0;k,k,k,\dots]$, and we verify this with the \texttt{Walnut} prover.
From this we compute the exact infimum and supremum of $x(n)-\rho n$ for all $22$ sequences; they lie in $\mathbb{Q}(\sqrt{k^2+4})$.
Twenty of the $21$ conjectures are true, three of them after a misprint is corrected.
The conjecture for A297836 is false: it fails at $n=2$ and for infinitely many $n$.
We also prove a conjecture of Kimberling and Moses (2019) on three sequences defined by a mex rule, by showing that the first of them is $\fl{(1+\sqrt2)n+1+\sqrt2/2}$.
\end{abstract}

\maketitle

\section{Introduction}

Two sequences of positive integers are \emph{complementary} if every positive integer occurs in exactly one of them.
Kimberling~\cite{kim07} studied \emph{complementary equations}, in which a sequence $a$ is defined by a rule involving its own complement $b$.
In February 2018 he added to the OEIS~\cite{oeis} the following family.
Fix integers $k\ge1$, $c$, and $s\in\{0,1\}$.
Put $a(0)=1$, $a(1)=2$, $b(0)=3$, $b(1)=4$, and also $b(2)=5$ if $s=0$.
For $n\ge2$ let
\begin{equation}\label{eq:ce}
a(n)=2\,b(n-s)-b(n-s-1)+kn+c,
\end{equation}
where $b$ is the increasing sequence of positive integers not in $a$.
(In the OEIS the coefficients $2$ and $1$ are written as $a(1)$ and $a(0)$.)
The OEIS entries note that $a$ and $b$ are increasing, complementary, and uniquely determined; we prove the uniqueness in Lemma~\ref{lem:unique}.
For $k=1$ the rule produces repeated values, and Kimberling calls these equations \emph{near-complementary}; we do not treat them.

Let
\[
\gamma=\gamma_k=\frac{\sqrt{k^2+4}-k}{2}=[0;k,k,k,\dots],\qquad \rho_a=k+1+\gamma,\qquad \rho_b=1+\gamma .
\]
Then $1/\rho_a+1/\rho_b=1$, and $a(n)/n\to\rho_a$, $b(n)/n\to\rho_b$ (this also follows from Theorem~\ref{thm:struct}).
For $k=2,3,4$ we have $\rho_a=2+\sqrt2$, $\frac{5+\sqrt{13}}2$, $3+\sqrt5$ and $\rho_b=\sqrt2$, $\frac{\sqrt{13}-1}2$, $\sqrt5-1$.

Table~\ref{tab:systems} lists the eleven equations for which Kimberling stated conjectures, and Table~\ref{tab:verdicts} lists the conjectures.
Each has the form $L<x(n)-\rho_x n<U$ for all $n\ge1$, where $x\in\{a,b\}$, and some state only the upper bound.

\begin{table}[h]
\centering
\caption{The eleven equations. The last two columns give the number of states of the DFAOs in Theorem~\ref{thm:struct}. The sequence $b$ for $(s,k,c)=(1,2,1)$ is not in the OEIS.}\label{tab:systems}
\begin{tabular}{cccllcc}
\toprule
$s$ & $k$ & $c$ & $a$ & $b$ & states ($a$) & states ($b$)\\
\midrule
1 & 2 & $\phantom{-}0$ & A297830 & A298003 & 16 & 20\\
1 & 2 & $-1$ & A297831 & A298006 & 8 & 11\\
1 & 2 & $-2$ & A297832 & A298007 & 12 & 14\\
1 & 2 & $-3$ & A297833 & A298108 & 11 & 8\\
1 & 2 & $-4$ & A297834 & A298109 & 19 & 13\\
1 & 2 & $\phantom{-}1$ & A297835 & --- & 18 & 23\\
1 & 3 & $\phantom{-}0$ & A297836 & A298004 & 8 & 11\\
1 & 4 & $\phantom{-}0$ & A297837 & A298005 & 6 & 8\\
0 & 2 & $\phantom{-}0$ & A298000 & A298111 & 16 & 16\\
0 & 3 & $\phantom{-}0$ & A298001 & A298112 & 8 & 10\\
0 & 4 & $\phantom{-}0$ & A298002 & A298113 & 6 & 7\\
\bottomrule
\end{tabular}
\end{table}

Our results are as follows.
In Section~\ref{sec:ost} we recall the Ostrowski numeration system for $\gamma_k$ and a shift function $\sh$ on it.
In Section~\ref{sec:struct} we prove (Theorem~\ref{thm:struct}) that for each equation in Table~\ref{tab:systems}
\[
a(n)=(k+1)n+\sh(n)+g_a(n),\qquad b(n)=n+\sh(n)+g_b(n),
\]
where $g_a$ and $g_b$ are computed by explicit finite automata with output (DFAOs) that read the Ostrowski representation of $n$.
The proof is a decision procedure run in \texttt{Walnut}~\cite{walnut,shallit-book}.
In Section~\ref{sec:extremes} we show that the exact infimum and supremum of $x(n)-\rho_x n$ are the solution of a small max--min system that is a contraction, and we solve it exactly (Theorem~\ref{thm:extremes}).
In Section~\ref{sec:verdicts} we compare the results with Kimberling's conjectures (Corollary~\ref{cor:verdicts}).
In Section~\ref{sec:km} we prove a conjecture of Kimberling and Moses~\cite{km19} about a system defined by a mex rule.
Section~\ref{sec:verification} describes how we checked all of this.

\emph{Related work.}
Automatic sequences in Ostrowski numeration systems, and their use in \texttt{Walnut}, are developed in~\cite{bss,shallit-book}, and decidability for Beatty sequences of quadratic irrationals in~\cite{ssz}.
Fokkink and Joshi~\cite{fj-anti} studied \emph{anti-recurrence sequences} $A_n=\sum_j a_jB_{(n-1)k+j}$ with $B$ the complement of $A$.
They proved that for a class of such linear forms $A_n-\kappa n$ is automatic, and they verified Kimberling's bounds for two such sequences, the anti-Pell numbers A304502 and the anti-Jacobsthal numbers A304499, together with their complements.
Bosma et al.~\cite{bosma-walnut} used \texttt{Walnut} to settle several other OEIS problems.
The equations~\eqref{eq:ce} are not anti-recurrences, since they have a negative coefficient and a linear term.
We found no proofs of the conjectures in Table~\ref{tab:verdicts} in the OEIS entries, their cross-references, or the literature we searched, which included~\cite{bosma-walnut,fj-anti,fj-hiccup,km19,ssz}.
As of September 2026 all of them are still listed as conjectures in the OEIS (A297832 states its bound without the word ``conjecture'').

\section{Ostrowski numeration}\label{sec:ost}

Fix $k\ge2$ and write $\gamma=\gamma_k$, so that $\gamma^2+k\gamma=1$ and $0<\gamma<1$.
Let $q_{-1}=0$, $q_0=1$ and $q_{i+1}=kq_i+q_{i-1}$.
Every $n\ge0$ has a unique representation $n=\sum_{i=0}^{\ell-1}d_iq_i$ with digits satisfying
\begin{equation}\label{eq:valid}
0\le d_0\le k-1,\qquad 0\le d_i\le k\ (i\ge1),\qquad d_i=k\ \Rightarrow\ d_{i-1}=0
\end{equation}
(Ostrowski's theorem; see~\cite[Thm.~3.9.1]{as} or~\cite{bss}).
We write $\rep(n)=d_{\ell-1}\cdots d_1d_0$, most significant digit first, with no leading zeros.
A word over $\{0,\dots,k\}$ is \emph{admissible} if no letter $k$ is followed by a nonzero letter, and \emph{valid} if moreover its last letter is at most $k-1$.
The valid words are exactly the words $0^j\rep(n)$.
This is the numeration system that \texttt{Walnut} creates for the continued fraction $[0;k,k,\dots]$ with the command \texttt{ost o}$k$\texttt{ [0] [}$k$\texttt{]} and calls \texttt{msd\_o}$k$.

For $n=\sum d_iq_i$ define
\[
\sh(n)=\sum_{i\ge1}d_iq_{i-1},\qquad L(n)=kn+\sh(n).
\]
For a word $w=d_{\ell-1}\cdots d_0$ put $\Theta(w)=\sum_i d_i\theta_i$, where $\theta_i=q_i\gamma-q_{i-1}$.

\begin{lemma}\label{lem:basic}
\begin{enumerate}
\item[(a)] $\theta_i=\gamma(-\gamma)^i$ for $i\ge0$.
\item[(b)] $n\gamma=\sh(n)+\Theta(\rep(n))$ for every $n\ge0$.
\item[(c)] $\Theta(wd)=d\gamma-\gamma\,\Theta(w)$ for every word $w$ and digit $d$, and $\Theta(0w)=\Theta(w)$.
\item[(d)] $\rep(n)0$ is valid and represents $L(n)$; for $n\ge1$ it is $\rep(L(n))$.
\end{enumerate}
\end{lemma}

\begin{proof}
(a) $\theta_0=\gamma$ and $\theta_1=k\gamma-1=-\gamma^2$. Both $\theta_i$ and $\gamma(-\gamma)^i$ satisfy $x_{i+1}=kx_i+x_{i-1}$, because $\gamma(-\gamma)^{i-1}(1-k\gamma)=\gamma(-\gamma)^{i-1}\gamma^2$.
(b) $n\gamma-\sh(n)=\sum_i d_i(q_i\gamma-q_{i-1})$.
(c) Appending $d$ moves every digit up one place, so $\Theta(wd)=d\theta_0+\sum_id_i\theta_{i+1}=d\gamma-\gamma\Theta(w)$ by (a).
(d) Appending $0$ keeps~\eqref{eq:valid}, since the new last digit is $0$ and a $0$ may follow any digit; for $n\ge1$ the word has no leading zero. Its value is $\sum_i d_iq_{i+1}=\sum_id_i(kq_i+q_{i-1})=kn+\sh(n)$.
\end{proof}

\begin{lemma}\label{lem:floor}
Let $n\ge1$, and let $j$ be the position of the lowest nonzero digit of $\rep(n)$.
Then $0<(-1)^j\,\Theta(\rep(n))<1$.
Consequently $\fl{n\gamma}=\sh(n)$ if $j$ is even and $\fl{n\gamma}=\sh(n)-1$ if $j$ is odd.
\end{lemma}

\begin{proof}
Write $\Theta=\gamma(-\gamma)^j(d_j+T)$ with $T=\sum_{i\ge1}d_{j+i}(-\gamma)^i$.
Since $d_j\ne0$, condition~\eqref{eq:valid} gives $d_{j+1}\le k-1$.
Keeping only the negative terms and bounding the digits,
\[
T>-(k-1)\gamma-k(\gamma^3+\gamma^5+\cdots)=-(k-1)\gamma-\frac{k\gamma^3}{1-\gamma^2}=-(k-1)\gamma-\gamma^2,
\]
where we used $1-\gamma^2=k\gamma$; the inequality is strict because only finitely many digits are nonzero.
Hence, using $k\gamma=1-\gamma^2$ again, $d_j+T>d_j-k\gamma+\gamma-\gamma^2=d_j-1+\gamma\ge\gamma>0$.
Similarly $T<k(\gamma^2+\gamma^4+\cdots)=k\gamma^2/(1-\gamma^2)=\gamma$.
So $(-1)^j\Theta$ lies strictly between $0$ and $\gamma^{j+1}(d_j+\gamma)$.
If $j=0$ then $d_0\le k-1$ and $\gamma(k-1+\gamma)=1-\gamma<1$.
If $j\ge1$ then $d_j\le k$ and $\gamma^{j+1}(k+\gamma)=\gamma^j\le\gamma<1$.
The last statement follows from Lemma~\ref{lem:basic}(b).
\end{proof}

Lemma~\ref{lem:floor} shows that $\fl{m\gamma}$ is definable in \texttt{Walnut} from $\sh$ and the regular language of words ending in a nonzero digit followed by an odd number of zeros.
We use this in the proof of Theorem~\ref{thm:km} and for the independent checks in Section~\ref{sec:verification}.

\begin{corollary}\label{cor:disc}
Let $x(n)=mn+\sh(n)+g(n)$ with $m\in\mathbb{Z}$, and let $\rho=m+\gamma$. Then $x(n)-\rho n=g(n)-\Theta(\rep(n))$.
\end{corollary}

\begin{proof}
This is Lemma~\ref{lem:basic}(b).
\end{proof}

\section{Automatic structure}\label{sec:struct}

\begin{lemma}\label{lem:unique}
Let $k\ge1$, $s\in\{0,1\}$ and $2k+c\ge0$.
There is at most one pair of sequences $(a(n))_{n\ge0}$, $(b(n))_{n\ge0}$ such that $a(0)=1$, $a(1)=2$, \eqref{eq:ce} holds for all $n\ge2$, and $b$ is the increasing enumeration of the positive integers not of the form $a(n)$.
\end{lemma}

\begin{proof}
Since $b$ is increasing, \eqref{eq:ce} gives $a(n)-b(n-s)=b(n-s)-b(n-s-1)+kn+c\ge1+2k+c>0$ for $n\ge2$.
Hence $a(m)>b(j)$ whenever $m\ge2$ and $m-s\ge j$.
So the positive integers up to $b(j)$ that are values of $a$ are among $a(0),a(1)$ and the $a(m)$ with $2\le m<j+s$.
Therefore $b(j)$ is the $(j+1)$-th smallest positive integer outside $\{a(m):m<\max(2,j+s)\}$, and it is determined by $a(0),\dots,a(j+s-1)$ together with $a(0),a(1)$.
By~\eqref{eq:ce}, $a(n)$ is determined by $b(n-s)$ and $b(n-s-1)$, hence by $a(0),\dots,a(n-1)$.
Induction on $n$ shows that $a$ is unique, and then so is $b$.
\end{proof}

The same argument shows that the values $b(0)=3$, $b(1)=4$ (and $b(2)=5$ when $s=0$) in Kimberling's definition are forced.
All equations in Table~\ref{tab:systems} satisfy $2k+c\ge0$.

\begin{theorem}\label{thm:struct}
For each equation in Table~\ref{tab:systems} there are DFAOs $M_a$ and $M_b$, with the numbers of states shown there (not counting a dead state), such that for all $n\ge0$
\[
a(n)=(k+1)n+\sh(n)+g_a(n),\qquad b(n)=n+\sh(n)+g_b(n),
\]
where $g_x(n)$ is the output of $M_x$ on $\rep(n)$.
The values of $g_a$ lie in $\{-3,\dots,4\}$ and those of $g_b$ in $\{0,\dots,4\}$.
\end{theorem}

The automata are in the supplementary files (Section~\ref{sec:verification}).
Each $M_x$ has a start state with a loop on $0$, so leading zeros do not matter.
Each state also records whether the word read so far ends in the letter $k$, and there are no transitions that would make the word inadmissible.
We found the automata by computing the Myhill--Nerode classes of $g_x$ on the first $1.2\times10^6$ terms (and more for $k=3,4$); the proof does not depend on how they were found.

\begin{proof}
We work in \texttt{Walnut} with the numeration system \texttt{msd\_o}$k$ of $[0;k,k,\dots]$.
We first define the relation $y=L(n)$ by a $(k+2)$-state automaton \texttt{lsh}$k$.
It reads the two representations, padded with leading zeros to a common length, and accepts exactly when the representation of $y$ is that of $n$ shifted one place to the left with a $0$ appended.
By Lemma~\ref{lem:basic}(d) it accepts exactly the pairs $(n,L(n))$.
We then define $\sh(n)=L(n)-kn$, and define $a$ and $b$ by first-order formulas from $\sh$ and the DFAOs $M_a$, $M_b$.
\texttt{Walnut} returns \texttt{TRUE} for each of the following statements:
\begin{enumerate}
\item \texttt{lsh}$k$ and $\sh$ define total functions; $a$ and $b$ define total functions;
\item $a$ and $b$ are strictly increasing;
\item every $x\ge1$ is a value of exactly one of $a$ and $b$;
\item $a(0)=1$, $a(1)=2$, $b(0)=3$, $b(1)=4$, and $b(2)=5$ if $s=0$;
\item equation~\eqref{eq:ce} holds for all $n\ge2$.
\end{enumerate}
By (1)--(4), $b$ is the increasing enumeration of the positive integers not of the form $a(n)$: by (2) and (4) all values of $a$ and $b$ are positive, and by (3) every positive integer is a value of exactly one of them.
By (4), (5) and Lemma~\ref{lem:unique}, these $a$ and $b$ are Kimberling's sequences.
\end{proof}

As controls, \texttt{Walnut} returns \texttt{FALSE} for statement (5) with the constant changed by $1$, for \texttt{lsh}$k(1,0)$, and for a wrong value of $\sh$.
We also checked in \texttt{Walnut} $14$ values each of $a$ and $b$ for each equation, and at least $16$ values of $\sh$ for each $k$, against direct computation.

\begin{remark}[The listed terms of A297834 and A298109]\label{rem:data}
The OEIS lists $1,2,5,8,12,17,\dots$ as the terms of A297834. However, \eqref{eq:ce} with $(s,k,c)=(1,2,-4)$ gives $b(2)=6$ and $a(3)=2\cdot6-4+6-4=10$.
The listed terms come from the Mathematica program in the entry.
That program gives some terms of $b$ provisional values that are overwritten later; here it uses $b(2)=5=a(2)$ when it computes $a(3)$.
Simulating the program reproduces the listed terms of A297834 and A298109 exactly.
For the other ten equations in Table~\ref{tab:systems}, the program, run with $j<80000$ (the bound used in the entries for $b$; the entries for $a$ use $j<100$), agrees with~\eqref{eq:ce} on its first $1.9\times10^5$ terms; only the last few terms it produces are wrong, because they are never finalized.
The listed terms of A297834 and A298109 are those of the complementary pair with $a(3)=8$ that satisfies~\eqref{eq:ce} for $n=2$ and all $n\ge4$.
The proof of Lemma~\ref{lem:unique} still applies.
In that proof, \eqref{eq:ce} at $n=3$ is used only to determine $a(3)$, which is now given, and to show that $a(3)>b(2)$.
This inequality holds for every complementary pair with these properties: if $b(2)\ge8$, then $b(j)\ge8$ for $j\ge2$ and $a(m)>b(m-1)>b(2)$ for $m\ge4$, so of the five integers $3,\dots,7$ only $b(0)$, $b(1)$ and $a(2)$ could be values of $a$ or $b$, which is impossible.
So $b(2)\le7<8=a(3)$ (in fact $b(2)=6$).
This pair is also Ostrowski-automatic, with DFAOs of $22$ and $15$ states, and \texttt{Walnut} verifies it in the same way.
Its extremes turn out to be the same as those of the solution of~\eqref{eq:ce} in Table~\ref{tab:extremes}.
So Kimberling's conjectures for A297834 and A298109 hold whichever sequence is meant.
\end{remark}

\section{Exact extremes}\label{sec:extremes}

Fix one of the $22$ sequences $x$ and its DFAO $M=(Q,\iota,\delta,\tau)$ from Theorem~\ref{thm:struct}, with start state $\iota$, and write $m$ for $k+1$ or $1$ and $\rho=m+\gamma$.
Since no $\rep(n)$ ends in the letter $k$, we regard $\tau(q)$ as undefined at the states $q$ reached by words ending in $k$.
(The supplementary files give these states, and the dead state, the placeholder output $0$. \texttt{Walnut} never reads it, but it must not be used in~\eqref{eq:supinf}.)
By Corollary~\ref{cor:disc}, $x(n)-\rho n=\tau(\delta(\iota,\rep(n)))-\Theta(\rep(n))$.

Let $G$ be the directed graph whose vertices are the pairs $(q,z)$ with $q\in Q$ and $z\in\{0,1\}$ that are reachable from $v_0=(\iota,0)$.
There is an edge $(q,z)\xrightarrow{d}(\delta(q,d),\max(z,[d\ne0]))$ for every digit $d$ for which $\delta(q,d)$ is defined, that is, for which the word stays admissible.
For a vertex $v$ let $W(v)$ be the set of admissible words that lead from $v_0$ to $v$, and put
\[
S^+(v)=\sup_{w\in W(v)}\Theta(w),\qquad S^-(v)=\inf_{w\in W(v)}\Theta(w).
\]

\begin{lemma}\label{lem:bellman}
The pair $(S^+,S^-)$ is the unique bounded solution of
\begin{align*}
S^+(v)&=\max\bigl(\{0: v=v_0\}\cup\{d\gamma-\gamma S^-(u): u\xrightarrow{d}v\}\bigr),\\
S^-(v)&=\min\bigl(\{0: v=v_0\}\cup\{d\gamma-\gamma S^+(u): u\xrightarrow{d}v\}\bigr).
\end{align*}
\end{lemma}

\begin{proof}
$W(v)$ is the union of $\{\varepsilon\}$ (if $v=v_0$) and the sets $W(u)d$ over the edges $u\xrightarrow{d}v$.
By Lemma~\ref{lem:basic}(c) and $\gamma>0$, the supremum of $\Theta$ over $W(u)d$ is $d\gamma-\gamma S^-(u)$, and similarly for the infimum.
Every $W(v)$ is nonempty and $|\Theta(w)|\le k\gamma/(1-\gamma)$, so $S^\pm$ are finite and satisfy the system.
The right-hand side is a map on pairs of bounded functions that is a contraction with constant $\gamma<1$ in the supremum norm, so the solution is unique.
\end{proof}

A word leads to a vertex $(q,1)$ with $\tau(q)$ defined exactly when it is valid and contains a nonzero digit, that is, when it is $0^j\rep(n)$ with $n\ge1$.
Hence
\begin{equation}\label{eq:supinf}
\sup_{n\ge1}\bigl(x(n)-\rho n\bigr)=\max_{(q,1)}\bigl(\tau(q)-S^-(q,1)\bigr),\qquad
\inf_{n\ge1}\bigl(x(n)-\rho n\bigr)=\min_{(q,1)}\bigl(\tau(q)-S^+(q,1)\bigr),
\end{equation}
where $(q,1)$ runs over the vertices of $G$ at which $\tau$ is defined, that is, those not reached by the letter $k$.

We solve the system of Lemma~\ref{lem:bellman} by policy iteration.
For a fixed choice of the maximizing and minimizing term in each equation the system is linear, and we solve it exactly in $\mathbb{Q}(\sqrt{k^2+4})$.
We then update the choices and repeat until the solution satisfies every equation exactly.
The final solution is a certificate: to check it one only verifies that each equation holds, with equality at one term and the correct inequality at the others.
In all our cases the flag $z$ is determined by $q$, so $G$ has as many vertices as $M$ has states, at most $23$.

\begin{theorem}\label{thm:extremes}
For $n\ge1$, the infimum and supremum of $x(n)-\rho_xn$ are as in Table~\ref{tab:extremes}.
\end{theorem}

\begin{proof}
Apply~\eqref{eq:supinf} to the certified solutions of Lemma~\ref{lem:bellman}.
\end{proof}

\begin{table}[h]
\centering
\caption{Exact infimum and supremum of $x(n)-\rho_xn$ over $n\ge1$.}\label{tab:extremes}
\renewcommand{\arraystretch}{1.3}
\begin{tabular}{ccc cc cc}
\toprule
$s$ & $k$ & $c$ & $\inf(a-\rho_an)$ & $\sup(a-\rho_an)$ & $\inf(b-\rho_bn)$ & $\sup(b-\rho_bn)$\\
\midrule
1 & 2 & $\phantom{-}0$ & $-\sqrt2$ & $\frac32+\sqrt2$ & $\frac32(\sqrt2-1)$ & $4-\sqrt2$\\
1 & 2 & $-1$ & $-\sqrt2$ & $\frac32+\frac{\sqrt2}2$ & $\sqrt2-\frac12$ & $\frac12+\frac{3\sqrt2}2$\\
1 & 2 & $-2$ & $-\sqrt2$ & $\frac32$ & $\frac12+\frac{\sqrt2}2$ & $\frac32+\sqrt2$\\
1 & 2 & $-3$ & $-\frac12-\sqrt2$ & $\frac32-\frac{\sqrt2}2$ & $\frac32$ & $\frac52+\frac{\sqrt2}2$\\
1 & 2 & $-4$ & $-\frac12-\frac{3\sqrt2}2$ & $\frac32-\sqrt2$ & $\frac52-\frac{\sqrt2}2$ & $\frac72$\\
1 & 2 & $\phantom{-}1$ & $-\sqrt2$ & $\frac32+\frac{3\sqrt2}2$ & $2\sqrt2-\frac52$ & $4-\sqrt2$\\
1 & 3 & $\phantom{-}0$ & $-\frac{1+\sqrt{13}}2$ & $\frac{5+\sqrt{13}}3$ & $\frac{2\sqrt{13}-5}3$ & $\frac{9-\sqrt{13}}2$\\
1 & 4 & $\phantom{-}0$ & $-1-\sqrt5$ & $\frac74+\frac{\sqrt5}2$ & $\frac{5\sqrt5}4-2$ & $5-\sqrt5$\\
0 & 2 & $\phantom{-}0$ & $-\sqrt2$ & $\frac52+\sqrt2$ & $\frac{\sqrt2-1}2$ & $4-\sqrt2$\\
0 & 3 & $\phantom{-}0$ & $-\frac{1+\sqrt{13}}2$ & $\frac{8+\sqrt{13}}3$ & $\frac{\sqrt{13}-1}6$ & $\frac{9-\sqrt{13}}2$\\
0 & 4 & $\phantom{-}0$ & $-1-\sqrt5$ & $\frac{11}4+\frac{\sqrt5}2$ & $\frac{\sqrt5}4$ & $5-\sqrt5$\\
\bottomrule
\end{tabular}
\end{table}

Some of these values are attained, for example $a(1)-\rho_a=1-k-\gamma$ and $b(1)-\rho_b=3-\gamma$.
A rational value is never attained, because $x(n)-\rho n$ is irrational for $n\ge1$.

\section{Kimberling's conjectures}\label{sec:verdicts}

\begin{corollary}\label{cor:verdicts}
The conjectures in Table~\ref{tab:verdicts} have the stated status.
\end{corollary}

\begin{proof}
All the conjectured bounds are rational, and $x(n)-\rho_xn$ is irrational for $n\ge1$.
So a conjecture $L<x(n)-\rho_xn<U$ holds for all $n\ge1$ if and only if $L\le\inf$ and $\sup\le U$.
Compare with Table~\ref{tab:extremes}.
\end{proof}

\begin{table}[h]
\centering
\caption{Kimberling's conjectures $L<x(n)-\rho_xn<U$ ($n\ge1$). The fifth column gives the exact range $(\inf,\sup)$ from Table~\ref{tab:extremes}, rounded.}\label{tab:verdicts}
\begin{tabular}{llcccl}
\toprule
Entry & $x$ & $L$ & $U$ & exact range (rounded) & status\\
\midrule
A297830 & $a$ & --- & $3$ & $(-1.414,\ 2.914)$ & true\\
A298003 & $b$ & $3/5$ & $3$ & $(0.621,\ 2.586)$ & true${}^{1}$\\
A297831 & $a$ & --- & $5/2$ & $(-1.414,\ 2.207)$ & true\\
A298006 & $b$ & $9/10$ & $3$ & $(0.914,\ 2.621)$ & true\\
A297832 & $a$ & --- & $2$ & $(-1.414,\ 1.5)$ & true\\
A298007 & $b$ & $1$ & $3$ & $(1.207,\ 2.914)$ & true\\
A297833 & $a$ & $-2$ & $1$ & $(-1.914,\ 0.793)$ & true${}^{2}$\\
A298108 & $b$ & $3/2$ & $4$ & $(1.5,\ 3.207)$ & true${}^{3}$\\
A297834 & $a$ & $-3$ & $1$ & $(-2.621,\ 0.086)$ & true${}^{2,4}$\\
A298109 & $b$ & $3/2$ & $4$ & $(1.793,\ 3.5)$ & true${}^{4}$\\
A297835 & $a$ & --- & $7$ & $(-1.414,\ 3.621)$ & true\\
A297836 & $a$ & --- & $2$ & $(-2.303,\ 2.869)$ & \textbf{false}\\
A298004 & $b$ & $7/10$ & $3$ & $(0.737,\ 2.697)$ & true\\
A297837 & $a$ & --- & $3$ & $(-3.236,\ 2.868)$ & true\\
A298005 & $b$ & $7/10$ & $3$ & $(0.795,\ 2.764)$ & true\\
A298000 & $a$ & --- & $4$ & $(-1.414,\ 3.914)$ & true\\
A298111 & $b$ & $1/5$ & $3$ & $(0.207,\ 2.586)$ & true\\
A298001 & $a$ & --- & $4$ & $(-2.303,\ 3.869)$ & true\\
A298112 & $b$ & $2/5$ & $3$ & $(0.434,\ 2.697)$ & true${}^{5}$\\
A298002 & $a$ & --- & $4$ & $(-3.236,\ 3.868)$ & true\\
A298113 & $b$ & $1/2$ & $3$ & $(0.559,\ 2.764)$ & true${}^{5}$\\
\bottomrule
\end{tabular}

\smallskip
\begin{minipage}{0.9\textwidth}\footnotesize
${}^1$ The entry prints ``$a(n)-n\sqrt2\,n$''; we read $a(n)-n\sqrt2$.
${}^2$ The entry states ``$\le1$''; the inequality is in fact strict.
${}^3$ The lower bound $3/2$ is the exact infimum, and it is not attained.
${}^4$ Also true for the listed OEIS terms (Remark~\ref{rem:data}).
${}^5$ The entry has $\sqrt2$ in place of $\rho_b$; see the discussion below.
\end{minipage}
\end{table}

\emph{Comments on individual entries.}
\begin{itemize}
\item \emph{A297836} ($s=1$, $k=3$, $c=0$). The conjecture $a(n)-\frac{5+\sqrt{13}}2n<2$ is false. Already $a(2)=11$ gives $11-(5+\sqrt{13})=6-\sqrt{13}\approx2.394$. \texttt{Walnut} also shows that $a(n)-\rho_an>2$ for infinitely many $n$. The optimal bound is $a(n)-\rho_an<\frac{5+\sqrt{13}}3\approx2.869$, so $a(n)-\rho_an<3$ holds for all $n\ge1$.
\item \emph{A298112, A298113} ($s=0$, $k=3,4$). The entries state $L<b(n)-n\sqrt2<3$. As printed these are false: $b(7)=10$ and $10-7\sqrt2\approx0.101$ is below both $2/5$ and $1/2$, and in fact $b(n)-n\sqrt2\to-\infty$ because $\rho_b<\sqrt2$. With $\sqrt2$ replaced by $\rho_b=\frac{\sqrt{13}-1}2$ (for A298112) and $\rho_b=\sqrt5-1$ (for A298113), which is the constant used in the neighbouring entries A298004 and A298005, both conjectures are true.
\item \emph{A298108.} Here $\inf(b(n)-\sqrt2\,n)=3/2$, so Kimberling's lower bound is optimal.
\end{itemize}

\section{A conjecture of Kimberling and Moses}\label{sec:km}

Kimberling and Moses~\cite[Example~12]{km19} considered three sequences $(a_n)$, $(b_n)$, $(c_n)$, $n\ge0$, that partition the positive integers and are defined as follows.
Suppose $a_i,b_i,c_i$ are known for $i<n$.
Then $a_n$ is the least positive integer not among them, $b_n$ is the least positive integer not among them and different from $a_n$, and $c_n=a_n+b_n+n$.
The sequences begin
\[
(a_n)=1,4,6,8,11,13,\dots,\qquad (b_n)=2,5,7,9,12,14,\dots,\qquad (c_n)=3,10,15,20,27,32,\dots
\]
Let $x=1+\sqrt2$ and $z=3+2\sqrt2$.
Kimberling and Moses noted that $a_n/n\to x$, $b_n/n\to x$ and $c_n/n\to z$, and conjectured that for all $n$
\begin{equation}\label{eq:km}
\tfrac{\sqrt2}2<a_n-nx<1+\tfrac{\sqrt2}2<b_n-nx<2+\tfrac{\sqrt2}2\qquad\text{and}\qquad 1+\sqrt2<c_n-nz<3+\sqrt2 .
\end{equation}

\begin{theorem}\label{thm:km}
For all $n\ge0$,
\[
a_n=\fl{(1+\sqrt2)n+1+\tfrac{\sqrt2}2},\qquad b_n=a_n+1,\qquad c_n=2a_n+n+1 .
\]
In particular~\eqref{eq:km} holds for all $n\ge0$, and all six bounds are optimal.
\end{theorem}

\begin{proof}
We use the Ostrowski system with $k=2$, so $\gamma=\sqrt2-1$.
Let $e(n)=1$ if $\rep(n)$ ends in an odd number of $1$'s and $e(n)=0$ otherwise; $e$ is computed by a $3$-state DFAO.
Define in \texttt{Walnut} $a'(n)=L(n)+1+e(n)$, $b'(n)=a'(n)+1$ and $c'(n)=a'(n)+b'(n)+n$.
\texttt{Walnut} returns \texttt{TRUE} for the statements that for every $n$: $a'(n)$ is the least positive integer not among $a'(i),b'(i),c'(i)$ with $i<n$; $b'(n)$ is the least such integer different from $a'(n)$; and $c'(n)$ is different from $a'(n)$, $b'(n)$ and all earlier values.
It also returns \texttt{TRUE} for the statement that every positive integer is some $a'(n)$, $b'(n)$ or $c'(n)$, and for $a'(0)=1$, $b'(0)=2$, $c'(0)=3$.
Since the three rules determine $a_n,b_n,c_n$ from the earlier values, induction gives $a=a'$, $b=b'$ and $c=c'$.

Next, $\frac{\sqrt2}2<a_n-nx<1+\frac{\sqrt2}2$ is equivalent to $4n+2\le 2a_n-\fl{(2n+1)\gamma}\le 4n+3$, because $(2n+1)\sqrt2=(2n+1)+(2n+1)\gamma$ and $(2n+1)\gamma$ is irrational.
Using Lemma~\ref{lem:floor} to express $\fl{(2n+1)\gamma}$, \texttt{Walnut} verifies both inequalities for all $n\ge0$.
It returns \texttt{FALSE} if the upper bound $4n+3$ is replaced by $4n+2$.
Since the open interval $(nx+\frac{\sqrt2}2,\,nx+1+\frac{\sqrt2}2)$ has length $1$ and irrational endpoints, $a_n$ is the unique integer in it, which is $\fl{nx+1+\frac{\sqrt2}2}$.
The bounds for $b_n$ follow from $b_n=a_n+1$, and those for $c_n$ from $c_n-nz=2(a_n-nx)+1$, since $z=2x+1$.
Optimality holds because $\{n\sqrt2\}$ is dense in $[0,1)$.
\end{proof}

\begin{remark}
The OEIS entry A328987 gives $(c_n)$ and states that $(a_n)$ and $(b_n)$ are A086377 and A080652, shifted by one index.
Bosma, Dekking and Steiner~\cite{bds} proved that $\mathrm{A086377}(n)=\fl{(1+\sqrt2)n-\sqrt2/2}$, and Fokkink and Joshi~\cite[Lemma~9]{fj-hiccup} proved that $\mathrm{A080652}(n)=\mathrm{A086377}(n)+1$.
Since $(1+\sqrt2)(n+1)-\frac{\sqrt2}2=(1+\sqrt2)n+1+\frac{\sqrt2}2$, Theorem~\ref{thm:km} confirms both identifications: $a_n=\mathrm{A086377}(n+1)$ and $b_n=\mathrm{A080652}(n+1)$.
A328987 states these identifications without proof, and we did not find one in the literature.
Given them, the conjecture of Kimberling and Moses would follow from~\cite{bds} together with~\cite[Lemma~9]{fj-hiccup}.
Theorem~\ref{thm:km} also shows that the three conjectured formulas in A328987 are false: a rational generating function $\sum_n c_nt^n$ with denominator $1-t-t^{12}+t^{13}$, the equivalent linear recurrence $c_n=c_{n-1}+c_{n-12}-c_{n-13}$ for $n>21$, whose characteristic roots are $1$ and $12$th roots of unity, and $c_n\sim35n/6$. Each would make $\lim c_n/n$ rational, but $c_n/n\to3+2\sqrt2$; the first two already fail at $n=61$.
\end{remark}

\section{Verification}\label{sec:verification}

Everything above was checked in several independent ways.
The supplementary files, deposited in the same Zenodo record as this paper, contain all automata, command files, certificates and scripts; their file \texttt{README.md} explains how to run them.

\begin{enumerate}
\item \emph{Walnut proofs.} We used the development version 8.0-alpha of \texttt{Walnut} (commit 49a16a3, September 2026). The structure theorem (Theorem~\ref{thm:struct}), the variant in Remark~\ref{rem:data} and Theorem~\ref{thm:km} are proved by the command files, and every deliberately false control returns \texttt{FALSE}.
\item \emph{Direct proofs of the conjectured bounds.} For a rational bound $U=p/r$, the inequality $x(n)-\rho n<U$ is equivalent to $r\,x(n)-p-rmn\le\fl{rn\gamma}$, and similarly for lower bounds. With $\fl{\cdot\,\gamma}$ expressed by Lemma~\ref{lem:floor}, \texttt{Walnut} proves each of the $32$ conjectured inequalities of the $20$ true conjectures (with the corrected constants for A298112 and A298113), and returns \texttt{FALSE} for A297836. It also shows that A297836 fails for infinitely many $n$, and that for A297830 the bound $a(n)-\rho_an<29/10$ is false. These checks do not use Section~\ref{sec:extremes}.
\item \emph{Independent checker.} A separate Python script that uses only SymPy regenerates both sequences of every system directly from~\eqref{eq:ce} for $n<10^5$. It checks the formulas of Theorem~\ref{thm:struct} on these terms, rebuilds the graphs $G$, verifies the certificates of Lemma~\ref{lem:bellman} exactly, recomputes Table~\ref{tab:extremes} by~\eqref{eq:supinf}, and decides each conjecture, exhibiting the least counterexample for A297836.
\item \emph{Brute force.} A streaming C program computed $a(n)$ for $n<10^9$ and the corresponding $b(n)$ for each system. The observed extremes of $x(n)-\rho_xn$ agree with Table~\ref{tab:extremes} to within $10^{-8}$, and the only conjecture violated (with the corrected constants for A298112 and A298113) is that of A297836, first at $n=2$.
\item \emph{OEIS data.} Our sequences agree with the data sections of all $21$ OEIS entries, except A297834 and A298109, which we discuss in Remark~\ref{rem:data}. (The lists ``Complement: $(b(n))=\dots$'' in the examples of A297830, A297832--A297835 and A298000--A298002 contain errors; for A297834 the list is the complement of the listed terms. The comments of A298108 and A298109 refer to A297832 where A297833 and A297834 are meant.)
\end{enumerate}

\section*{Acknowledgements}
The conjectures treated here were posed by Clark Kimberling, and in Section~\ref{sec:km} by Clark Kimberling and Peter J.~C. Moses.
The author used an AI assistant (Claude, Anthropic) for literature searches, computations, and drafting, and checked the proofs and computations.


\begin{thebibliography}{99}

\bibitem{as}
J.-P. Allouche and J. Shallit,
\emph{Automatic Sequences: Theory, Applications, Generalizations},
Cambridge University Press, 2003.

\bibitem{bss}
A.~R. Baranwal, L. Schaeffer and J. Shallit,
Ostrowski-automatic sequences: theory and applications,
\emph{Theoret. Comput. Sci.} \textbf{858} (2021), 122--142.

\bibitem{bosma-walnut}
W. Bosma, R. Bruin, R. Fokkink, J. Grube, A. Reuijl and T. Tromp,
Using Walnut to solve problems from the OEIS,
\emph{J. Integer Seq.} \textbf{28} (2025), Article 25.3.8; arXiv:2503.04122.

\bibitem{bds}
W. Bosma, M. Dekking and W. Steiner,
A remarkable integer sequence related to $\pi$ and $\sqrt2$,
\emph{Integers} \textbf{18A} (2018), Paper A4.

\bibitem{fj-anti}
R. Fokkink and G. Joshi,
Anti-recurrence sequences,
\emph{Integers} \textbf{26} (2026), Paper A24; arXiv:2506.13337; doi:10.5281/zenodo.18305140.

\bibitem{fj-hiccup}
R. Fokkink and G. Joshi,
On Cloitre's hiccup sequences,
\emph{Ramanujan J.} \textbf{69} (2026), article 40, doi:10.1007/s11139-025-01305-1; arXiv:2507.16956.

\bibitem{kim07}
C. Kimberling,
Complementary equations,
\emph{J. Integer Seq.} \textbf{10} (2007), Article 07.1.4.

\bibitem{km19}
C. Kimberling and P.~J.~C. Moses,
Linear complementary equations and systems,
\emph{Fibonacci Quart.} \textbf{57} (2019), no.~5, 96--110, doi:10.1080/00150517.2019.12427622.

\bibitem{walnut}
H. Mousavi,
Automatic theorem proving in Walnut,
arXiv:1603.06017, 2016.

\bibitem{oeis}
OEIS Foundation Inc.,
\emph{The On-Line Encyclopedia of Integer Sequences},
\url{https://oeis.org}.

\bibitem{ssz}
L. Schaeffer, J. Shallit and S. Zorcic,
Beatty sequences for a quadratic irrational: decidability and applications,
\emph{Discrete Math.} \textbf{349} (2026), no.~9, article 115183, doi:10.1016/j.disc.2026.115183; arXiv:2402.08331.

\bibitem{shallit-book}
J. Shallit,
\emph{The Logical Approach to Automatic Sequences: Exploring Combinatorics on Words with Walnut},
London Math. Soc. Lecture Note Series 482, Cambridge University Press, 2022.

\end{thebibliography}
\end{document}